\documentclass[12pt]{amsart}

\usepackage[english]{babel}
\usepackage[utf8]{inputenc}
\usepackage[T1]{fontenc}
\usepackage{mathptmx} 
\usepackage{amsmath, amssymb, amsthm, amsfonts, amstext, amscd}
\usepackage{graphicx}
\usepackage[figuresright]{rotating}
\usepackage{longtable, booktabs}
\usepackage{caption}
\usepackage{multirow}
\usepackage{tensor}
\usepackage{enumerate}
\usepackage[table]{xcolor}
\usepackage{color}
\usepackage{fancyhdr}
\usepackage{hyperref}

\hypersetup{
	pdftitle={DC},
	pdfauthor={CC},
	pdfsubject={Eigenvalue},
	pdfkeywords={Laplacian, eigenvalues, Weyl's law, Riesz mean, universal estimates},
	pdfpagemode=UseOutlines,
	colorlinks=true,
	linkcolor=blue,
	urlcolor=cyan,
	citecolor=red
}

\newcommand{\norm}[1]{\left\|#1\right\|}

\newcommand{\ps}[2]{\langle{#1},{#2}\rangle}

\newcommand{\abs}[1]{\left|#1\right|}

\newcommand{\br}[1]{\left(#1\right)}

\newcommand{\dis}{\displaystyle}
\newcommand{\Lc}{\mathcal{L}}
\newcommand{\Lon}{\Lc_{0,n}^{cl}}
\newcommand{\Lrn}{\Lc_{\rho,n}^{cl}}

\def\dfrac{\displaystyle\frac}

\newtheorem{prop}{Proposition}[section]
\newtheorem{thm}{Theorem}[section]
\newtheorem{lem}{Lemma}[section]
\newtheorem{rem}{Remark}[section]
\newtheorem{cor}{Corollary}[section]

\newtheorem{conjecture}{Conjecture}[section]

\numberwithin{equation}{section}

\begin{document}
	
	\title{Estimates for Eigenvalues of the Dirichlet Laplacian on Riemannian Manifolds}
	\author[Daguang Chen, Qing-Ming Cheng]{Daguang Chen and Qing-Ming Cheng}
	\thanks{The first author is supported by NSFC grant No.11831005 and NSFC-FWO grant No.11961131001.\\
		\indent The second author was partially supported by JSPS Grant-in-Aid for Scientific Research: No. 25K06992.}
	\maketitle
	
	\begin{abstract}
		We revisit the eigenvalue problem of the Dirichlet Laplacian on bounded domains in complete Riemannian manifolds. 
		By building on classical results like Li-Yau's and Yang's inequalities, we derive upper and lower bounds for eigenvalues.  
		For the projective spaces and their minimal submanifolds, we also give explicit estimates on lower bounds for eigenvalues of the Dirichlet Laplacian.
	\end{abstract}
	
	\section{Introduction}
	
	Let \(\Omega\) be a bounded domain in an \(n\)-dimensional complete Riemannian manifold \(M\) with boundary (possibly empty). 
	The eigenvalue problem of the Dirichlet Laplacian on \(\Omega\) is given by
	\begin{equation}\label{DL}
		\begin{cases}
			\Delta u = -\lambda u, & \text{in } \Omega, \\
			u = 0, & \text{on } \partial \Omega,
		\end{cases}
	\end{equation}
	where \(\Delta\) denotes the Laplacian on \(M\). 
	It is well known that the spectrum of this problem consists of real, discrete eigenvalues
	\begin{equation*}
		0 < \lambda_1 < \lambda_2 \leq \lambda_3 \leq \cdots \nearrow \infty,	
	\end{equation*}
	where each eigenvalue \(\lambda_i\) has finite multiplicity and is repeated according to its multiplicity.
	
	The eigenvalue problem of the Dirichlet Laplacian arises from various problems of mathematical physics. It may refer to modes of an idealized drum, a mode of an idealized optical fiber in the paraxial approximation, as well as to small waves at the surface of an idealized pool. 
	
	It is well-known that for the eigenvalue problem of the Dirichlet Laplacian, we have Weyl's asymptotic formula \cite{Weyl11}  
	\begin{equation}\label{ineq:Weyl}
		\lambda_k \sim \frac{4\pi^2}{(\omega_n \abs{\Omega})^{\frac{2}{n}}} k^{\frac{2}{n}}, \quad k \to \infty,	
	\end{equation}  
	where \(\omega_n\) and \(\abs{\Omega}\) denote the volume of the unit ball in \(\mathbb{R}^n\) and \(\Omega\), respectively. 
	
	The paper is organized as follows. In Section \ref{section2}, we review eigenvalues of the Dirichlet Laplacian on a bounded domain in \(\mathbb{R}^n\). 
	In Section \ref{section3}, we revisit universal inequalities for eigenvalues of the Dirichlet Laplacian on a bounded domain in Riemannian manifolds and present the main theorems. In Section \ref{sec:pre}, we recall Karamata's Tauberian theorem and the first standard embedding of the projective spaces into Euclidean space, and provide the proofs for Theorem \ref{ChenCheng24} and Corollary \ref{cor:PS}.
	
	\section{Eigenvalues of Laplacian on a Bounded Domain in \(\mathbb{R}^n\)}\label{section2}
	
	When \(\Omega\) is a bounded domain in \(\mathbb{R}^n\), the study of universal inequalities for the eigenvalues of \eqref{DL} was initiated by Payne, P\'olya, and Weinberger in their seminal works \cite{PPW55,PPW56}. They established the following inequality
	\begin{equation}\label{ineq:PPW}
		\lambda_{k + 1} - \lambda_k \leq \frac{4}{nk} \sum_{i = 1}^{k} \lambda_i.	
	\end{equation}
	In 1980, Hile and Protter \cite{HP80} improved on the inequality of Payne, P\'olya, and Weinberger by proving
	\begin{equation}\label{ineq:HP80}
		\sum_{i = 1}^{k} \frac{\lambda_i}{\lambda_{k + 1} - \lambda_i} \geq \frac{nk}{4}.	
	\end{equation}
	It is very important to find a sharp universal inequality for eigenvalues in some sense. For this purpose, Yang \cite{Yang91} (cf. \cite{ChengYang07}) made a landmark work. He proved a significant inequality
	\begin{equation}\label{ineq:Yang1}
		\sum_{i = 1}^{k} (\lambda_{k + 1} - \lambda_i)^2 \leq \dfrac{4}{n} \sum_{i = 1}^{k} (\lambda_{k + 1} - \lambda_i) \lambda_i.	
	\end{equation}
	We should remark that the coefficient \(\frac{4}{n}\) is best possible according to Weyl's asymptotic formula (\ref{ineq:Weyl}), which cannot be improved.
	From \eqref{ineq:Yang1}, one can deduce
	\begin{equation}\label{ineq:Yang2}
		\lambda_{k + 1} \leq \frac{1}{k} \left(1 + \frac{4}{n}\right) \sum_{i = 1}^{k} \lambda_i.	
	\end{equation}
	These inequalities \eqref{ineq:Yang1} and \eqref{ineq:Yang2} are referred to as Yang's first and second inequalities, respectively (cf. \cite{Ashbaugh00,AshBen07,Ashbaugh17}). Using Chebyshev's inequality, it is straightforward to show the following logical relationships
	\begin{equation*}
		\eqref{ineq:Yang1} \implies \eqref{ineq:Yang2} \implies \eqref{ineq:HP80} \implies \eqref{ineq:PPW}.
	\end{equation*}
	
	The following is the famous P\'olya conjecture.
\vskip2mm
\noindent
{\bf P\'olya Conjecture}. 
		For the eigenvalue problem of the Dirichlet Laplacian, eigenvalues satisfy 
		\[
		\lambda_k \geq \frac{4\pi^2}{(\omega_n \abs{\Omega})^{\frac{2}{n}}} k^{\frac{2}{n}}, \quad \text{for } k = 1, 2, \dots.
		\] 
		P\'olya \cite{PolyaConj} resolved the case that the bounded domains are tiled. Furthermore, when \(\Omega\) is a ball, the P\'olya conjecture has been resolved by Filonov, Levitin, Polterovich, and Sher \cite{FLPS23}, very recently. In 1983, Li and Yau \cite{LY83} made significant progress toward resolving P\'olya's conjecture by making use of Fourier transformation. In the sense of summation, their result is best possible, that is, they proved 
	\begin{equation}\label{ineq:LiYau-sum}
		\frac{1}{k} \sum_{j=1}^k \lambda_j \geq \frac{n}{n+2} \frac{4\pi^2}{(\omega_n \abs{\Omega})^{\frac{2}{n}}} k^{\frac{2}{n}}, \quad \text{for } k = 1, 2, \dots,
	\end{equation}
	which implies 
	\begin{equation}\label{ineq:LiYau}
		\lambda_k \geq \frac{n}{n + 2} \frac{4\pi^2}{(\omega_n \abs{\Omega})^{\frac{2}{n}}} k^{\frac{2}{n}}, \quad \text{for } k = 1, 2, \dots.	
	\end{equation}
	According to the results of Li and Yau \cite{LY83} and the P\'olya conjecture, we know that the lower bounds for eigenvalue \(\lambda_k\) are given. As one sees, the lower bounds depend on the domain \(\Omega\). 
	
	On the other hand, the study on upper bounds for eigenvalues of the Laplacian was very difficult. In \cite{ChengYang07}, the second author and Yang successfully studied upper bounds for eigenvalues of the Laplacian. They obtained an upper bound of eigenvalue \(\lambda_{k+1}\) by the first eigenvalue \(\lambda_1\), and this inequality is universal, which does not depend on the domain.	
	\begin{thm}\label{thm:ChengYang07}
		For the eigenvalue problem of the Dirichlet Laplacian \eqref{DL}, eigenvalues satisfy
		\[
		\lambda_{k+1} \leq C_0(n, k) k^{\frac{2}{n}} \lambda_1,
		\]  
		where  
		\[
		C_0(n, k) =
		\begin{cases}
			\frac{j_{n/2,1}^2}{j_{n/2 - 1,1}^2}, & \text{for } k = 1, \\
			1 + \frac{a(\min\{n, k - 1\})}{n}, & \text{for } k \geq 2,
		\end{cases}
		\]  
		and \(a(1) \leq 2.64\) and \(a(m) \leq 2.2 - 4 \log\left(1 + \frac{m - 3}{50}\right)\) for \(m \geq 2\) are constants depending only on \(m\). Here, \(j_{p,k}\) denotes the \(k\)-th positive zero of the standard Bessel function \(J_p(x)\) of the first kind of order \(p\).
	\end{thm}
	\begin{rem} 
		From Weyl's asymptotic formula \eqref{ineq:Weyl}, it is clear that the upper bound obtained by the second author and Yang \cite{ChengYang07} is optimal in terms of the order of \(k\).  
	\end{rem}
	\begin{rem}
		In \cite{Ashbaugh17}, Professor S. Ashbaugh wrote that Cheng and Yang made great strides in the field, in what amounted to a tour de force in 2007.
	\end{rem}
	In order to prove their theorem \ref{thm:ChengYang07}, Cheng and Yang proved a recursion formula.
	\begin{thm}[The recursion formula of Cheng and Yang \cite{ChengYang07}]
		Let \(\mu_1 \leq \mu_2 \leq \dots \leq \mu_{k+1}\) be any positive real numbers satisfying 
		\begin{equation*}
			\sum_{i=1}^k (\mu_{k+1} - \mu_i)^2 \leq \frac{4}{t} \sum_{i=1}^k \mu_i (\mu_{k+1} - \mu_i).
		\end{equation*} 
		Define
		\begin{equation*}
			\begin{aligned}
				& G_k = \frac{1}{k} \sum_{i=1}^k \mu_i, \quad T_k = \frac{1}{k} \sum_{i=1}^k \mu_i^2, \\
				& F_k = \left(1 + \frac{2}{t}\right) G_k^2 - T_k.
			\end{aligned}
		\end{equation*}
		Then, we have, for any \(\ell, k\),
		\begin{equation}\label{eq.2-7}
			\frac{F_{k+\ell}}{(k+\ell)^{\frac{4}{t}}} \leq \frac{F_{k}}{k^{\frac{4}{t}}}.
		\end{equation} 
	\end{thm} 
	\begin{proof}[Proof of Theorem \ref{thm:ChengYang07}]
		According to Yang's first inequality (\ref{ineq:Yang1}), we know that eigenvalues \(\lambda_k\) satisfy the condition in the theorem with \(t = n\). By making use of the recursion formula of Cheng and Yang (\ref{eq.2-7}), we have
		\begin{equation}
			F_k \leq C(n,k-1) \left(\frac{k}{k-1}\right)^{\frac{4}{n}} F_{k-1} \leq k^{\frac{4}{n}} F_1 = \frac{2}{n} k^{\frac{4}{n}} \lambda_1^2.
		\end{equation}
		By making use of Yang's first inequality again, we obtain
		\[
		\left[\lambda_{k+1} - \left(1 + \frac{2}{n}\right) G_k\right]^2 \leq \left(1 + \frac{4}{n}\right) F_k - \frac{2}{n} \left(1 + \frac{2}{n}\right) G_k^2.
		\]
		Hence, we have
		\[
		\frac{\frac{2}{n}}{\left(1 + \frac{4}{n}\right)} \lambda_{k+1}^2 + \frac{1 + \frac{2}{n}}{1 + \frac{4}{n}} \left(\lambda_{k+1} - \left(1 + \frac{4}{n}\right) G_k\right)^2 \leq \left(1 + \frac{4}{n}\right) F_k.
		\]
		Thus, we derive
		\begin{equation}
			\lambda_{k+1}^2 \leq \frac{n}{2} \left(1 + \frac{4}{n}\right)^2 F_k \leq \left(1 + \frac{4}{n}\right)^2 k^{\frac{4}{n}} \lambda_1^2.
		\end{equation}
	\end{proof}
	
	For \(z \geq 0\), the Riesz mean of order \(\rho \ (\rho > 0)\), is defined as
	\begin{equation}\label{defn:RieszMean}
		R_{\rho}(z) = \sum_{k} \left(z - \lambda_k\right)_{+}^{\rho},	
	\end{equation}
	where \(\left(z - \lambda\right)_{+} := \max\left(0, z - \lambda\right)\) is the ramp function.
	As \(\rho \to 0^+\), the Riesz mean converges to the counting function
	\begin{equation}\label{defn:counting}
		N(z) = \sum_{\lambda_k \leq z} 1 = \sup_{\lambda_k \leq z} k.
	\end{equation}
	In terms of the counting function \eqref{defn:counting}, the Li-Yau inequality \eqref{ineq:LiYau-sum} states that
	\begin{equation}\label{li-yau-counting}
		N(z) \leq \left(\frac{n+2}{n}\right)^{n/2} \Lc_{0,n}^{cl} |\Omega| z^{n/2},
	\end{equation}
	where \(\Lc_{\rho,n}^{cl}\) is called the classical constant defined by
	\begin{equation}\label{defn:classical-constant}
		\Lrn = \dfrac{\Gamma(1 + \rho)}{\left(4 \pi\right)^{n/2} \Gamma(1 + \rho + n/2)}. 
	\end{equation}
	For \(\rho \geq 1\), Berezin \cite{Berezin73} proved that the Riesz means for the Dirichlet Laplacian satisfy
	\begin{equation}\label{ineq:laptev-weidl} 
		R_{\rho} (z) \leq \Lc_{\rho,n}^{cl} |\Omega| z^{\rho + n/2}.
	\end{equation}
	In \cite{LapWeidl1} (see also \cite{LapWeidl2}), Laptev and Weidl refer to \eqref{ineq:laptev-weidl} as the Berezin-Li-Yau inequality. In fact, they \cite{LapWeidl1} demonstrated the equivalence between the Li-Yau inequality \eqref{ineq:LiYau-sum} and the Berezin inequality \eqref{ineq:laptev-weidl} via the Legendre transform.
	
	Another well-known function associated with the spectrum is the trace of the heat kernel (equivalently, the partition function), denoted by \(Z(t)\). We recall the asymptotic formula of Kac \cite{Kac51} for \(Z(t)\)
	\begin{equation}\label{ineq:Kac51}
		Z(t) := \sum_{k = 1}^{\infty} e^{-\lambda_{k} t} \sim \frac{|\Omega|}{(4\pi t)^{n/2}},	
	\end{equation}
	which is equivalent to \eqref{ineq:Weyl} in terms of the Laplace transform. In \cite{Kac51}, Kac also established the inequality
	\begin{equation}\label{ineq:Kac}
		Z(t) = \sum_{k = 1}^{\infty} e^{-\lambda_{k} t} \leq \frac{|\Omega|}{(4\pi t)^{n/2}}.
	\end{equation}
	This result was refined in \cite{HS97}, where it was shown that \(t^{n/2} Z(t)\) is a nonincreasing function for \(t \to 0^+\). In \cite{HarHer11}, Harrell and Hermi showed that \eqref{ineq:laptev-weidl} is equivalent to \eqref{ineq:Kac} for \(\rho \geq 2\) via the Laplace transform.
	
	\section{Eigenvalues of Laplacian on a Bounded Domain in Riemannian Manifolds}\label{section3}
	
	Since Weyl's asymptotic formula \cite{Weyl11} for bounded domains in complete Riemannian manifolds holds, it is natural and important to derive universal inequalities for eigenvalues of the eigenvalue problem of the Dirichlet Laplacian on a bounded domain in a complete Riemannian manifold. 
	
	For the eigenvalue problem of the Dirichlet Laplacian on a compact homogeneous Riemannian manifold or on a compact minimal submanifold in a sphere, many mathematicians have studied universal inequalities for eigenvalues (see, for example, \cite{ChengYang05, ChengYang06, Cheng75, HS97, YangYau80} and others). Cheng and Yang \cite{ChengYang05, ChengYang06} derived optimal universal inequalities for eigenvalues of the Dirichlet Laplacian on a domain in a sphere or in a complex projective space. Namely, they proved
	\begin{thm}\label{thm:ChengYang05}
		For the eigenvalue problem of the Dirichlet Laplacian \eqref{DL} on a domain in the unit sphere, eigenvalues \(\lambda_k\) satisfy
		\begin{equation}
			\sum_{i=1}^{k} (\lambda_{k+1} - \lambda_i)^2 \leq \frac{4}{n} \sum_{i=1}^{k} (\lambda_{k+1} - \lambda_i) \br{\lambda_i + \frac{n^2}{4}}.
		\end{equation}
	\end{thm}
	\begin{rem}
		When \(\Omega \to \mathbb{S}^{n}(1)\), the above inequalities for all \(k\) become equalities. Hence, results of Cheng and Yang are optimal.
	\end{rem}
	Furthermore, since a sphere can be seen as a hypersurface in Euclidean space, Chen and Cheng \cite{ChenCheng08} studied the more general case of \(n\)-dimensional complete submanifolds in Euclidean space. They proved 
	\begin{thm}\label{thm:ChenCheng08submfd}
		Let \(\Omega\) be a bounded domain in an \(n\)-dimensional complete Riemannian manifold \(M^n\) isometrically immersed in the Euclidean space \(\mathbb{R}^N\). For the eigenvalue problem of the Dirichlet Laplacian \eqref{DL}, eigenvalues \(\lambda_k\) satisfy
		\begin{equation}\label{ineq:yang-first}
			\sum_{i=1}^{k} (\lambda_{k+1} - \lambda_i)^2 \leq \frac{4}{n} \sum_{i=1}^{k} (\lambda_{k+1} - \lambda_i) \br{\lambda_i + \frac{n^2}{4} H_0^2},
		\end{equation}
		where \(H\) is the mean curvature vector field of \(M^n\) with \(H_0^2 = \|H\|^2_{L^\infty(\Omega)} = \sup_{\Omega} |H|^2\).
	\end{thm}
	In order to prove our results, the following theorem of Cheng and Yang \cite{ChenCheng08} will play an important role.
	\begin{thm}\label{thm:chengyang06}
		Let \(\lambda_i\) be the \(i^{\text{th}}\) eigenvalue of the eigenvalue problem of the Dirichlet Laplacian on an \(n\)-dimensional compact Riemannian manifold \(\bar{\Omega} = \Omega \cup \partial \Omega\) with boundary \(\partial \Omega\) and \(u_i\) be the orthonormal eigenfunction corresponding to \(\lambda_i\). Then, for any function \(f \in C^3(\Omega) \cap C^2(\partial \Omega)\) and any integer \(k\), we have 
		\[
		\sum_{i=1}^k (\lambda_{k+1} - \lambda_i)^2 \|u_i \nabla f\|^2 \leq \sum_{i=1}^k (\lambda_{k+1} - \lambda_i) \norm{2 \nabla f \cdot \nabla u_i + u_i \Delta f}^2, 
		\]
		where \(\|f\|^2 = \dis \int_{M} f^2\) and \(\nabla f \cdot \nabla u_i = g(\nabla f, \nabla u_i)\).
	\end{thm}
	Let \(\Omega \subset M^n\) be a bounded domain and \(p \in \Omega\) be an arbitrary point of \(\Omega\) with a coordinate system \((x^1, \dots, x^n)\) in a neighborhood \(U\) of \(p\) in \(M^n\). Since \(M^n\) is an \(n\)-dimensional complete Riemannian manifold isometrically immersed in \(\mathbb{R}^N\), we can assume that \(y\) with components \(y^{\alpha}\) defined by 
	\[
	y^\alpha = y^\alpha (x^1, \dots, x^n), \qquad 1 \leq \alpha \leq N,
	\]
	is the position vector of \(p\) in \(\mathbb{R}^N\). We have 
	\[
	g_{ij} = g \br{\frac{\partial}{\partial x^i}, \frac{\partial}{\partial x^j}} = \ps{\sum_{\alpha=1}^N \frac{\partial y^\alpha}{\partial x^i} \frac{\partial}{\partial y^\alpha}}{\sum_{\beta=1}^N \frac{\partial y^\beta}{\partial x^j} \frac{\partial}{\partial y^\beta}} = \sum_{\alpha=1}^N \frac{\partial y^{\alpha}}{\partial x^i} \frac{\partial y^{\alpha}}{\partial x^j},
	\]
	where \(g\) denotes the induced metric of \(M^n\) from \(\mathbb{R}^N\), \(\langle \ , \ \rangle\) is the standard inner product in \(\mathbb{R}^N\). 
	\begin{lem}[{\cite[Lemma 2.1]{ChenCheng08}}]\label{lem:ChenCheng08}
		Let \(M\) be an \(n\)-dimensional complete Riemannian manifold with metric \(g\) isometrically immersed in a Euclidean space \(\mathbb{R}^N\). For any point \(p\) in \(M\), assuming that \(y\) with components \(y^{\alpha}\) defined by \(y^{\alpha} = y^{\alpha}(x^1, x^2, \dots, x^n)\) is the position vector of \(p\) in \(\mathbb{R}^N\), we have,
		\[
		\sum_{\alpha=1}^N g(\nabla y^{\alpha}, \nabla y^{\alpha}) = n, \quad \sum_{\alpha=1}^N (\Delta y^{\alpha})^2 = n^2 |H|^2,
		\]
		\[
		\sum_{\alpha=1}^N \Delta y^{\alpha} \nabla y^{\alpha} = 0, \quad \sum_{\alpha=1}^N g(\nabla y^{\alpha}, \nabla u)^2 = |\nabla u|^2,
		\]
		for any function \(u \in C^{\infty}(M)\), where \(H\) is the mean curvature vector of \(M\).
	\end{lem}
	\begin{proof}[Proof of Theorem \ref{thm:ChenCheng08submfd}]
		Let \(u_i\) be the eigenfunction corresponding to the eigenvalue \(\lambda_i\) such that \(\{u_i\}_{i \in \mathbb{N}}\) becomes an orthonormal basis of \(L^2(\Omega)\). Put \(f^{\alpha} = y^{\alpha}, 1 \leq \alpha \leq N\). Since \(M^n\) is complete and \(\Omega\) is a bounded domain, we know that \(\bar{\Omega}\) is a compact Riemannian manifold with boundary. From the theorem \ref{thm:chengyang06} of Cheng and Yang, we infer
		\[
		\sum (\lambda_{k+1} - \lambda_i)^2 \|u_i \nabla f^{\alpha}\|^2 \leq \sum (\lambda_{k+1} - \lambda_i) \|2 \nabla f^{\alpha} \cdot \nabla u_i + u_i \Delta f^{\alpha}\|^2. 
		\]
		Taking summation on \(\alpha\) and using the lemma \ref{lem:ChenCheng08}, we finish the proof.
	\end{proof}
	\begin{rem}
		Our results are optimal since for the unit sphere, \(H_0^2 = 1\), our inequality becomes one of Cheng and Yang.
	\end{rem}
	\begin{rem}
		Inequality \eqref{ineq:yang-first} had also been proved by El Soufi, Harrell, and Ilias \cite{SoufiHaIlias09}, independently.
	\end{rem}
	For a bounded domain in an \(n\)-dimensional complete Riemannian manifold isometrically minimally immersed in Euclidean space, we have 
	\begin{cor}[{\cite{ChenCheng08}}]\label{cor:CC08}
		Let \(\Omega\) be a bounded domain in an \(n\)-dimensional complete Riemannian manifold \(M^n\) isometrically minimally immersed in \(\mathbb{R}^N\). Then, for the eigenvalue problem (\ref{DL}), we have
		\begin{equation}\label{ineq:minimal}
			\sum_{i=1}^{k} (\lambda_{k+1} - \lambda_i)^2 \leq \frac{4}{n} \sum_{i=1}^k (\lambda_{k+1} - \lambda_i) \lambda_i,	
		\end{equation}
	\end{cor}
	\begin{rem}
		We would like to remark that Yang's first inequality does not only hold for domains in Euclidean spaces but also holds for domains in complete minimal submanifolds in Euclidean spaces.
	\end{rem}
	We should remark that the theorem \ref{thm:ChenCheng08submfd} of Chen and Cheng includes all complete Riemannian manifolds according to Nash's theorem. In fact, Cheng and Yang \cite{ChengYang11} obtained the following, by making use of a simple observation and the theorem \ref{thm:ChenCheng08submfd}.  
	\begin{thm}\label{ChengYang11}
		Let \(\Omega\) be a bounded domain in an \(n\)-dimensional complete Riemannian manifold \(M^n\). For the eigenvalue problem of the Dirichlet Laplacian (\ref{DL}), there exists a constant \(H_0^2\), which depends only on \(M\) and \(\Omega\), such that 
		\begin{equation}\label{CC08-1}
			\sum_{i=1}^{k} (\lambda_{k+1} - \lambda_i)^2 \leq \frac{4}{n} \sum_{i=1}^k (\lambda_{k+1} - \lambda_i) \left(\lambda_i + \frac{n^2}{4} H_0^2 \right).
		\end{equation}
	\end{thm}
	\begin{proof}
		According to Nash's theorem, we know that \(M^n\) can be immersed into Euclidean space \(\mathbb{R}^N\) by \(\varphi: M^n \to \mathbb{R}^N\). From the theorem \ref{thm:ChenCheng08submfd} of Chen and Cheng and putting \(H_0^2 = \dis \inf_{\varphi} \|H\|^2_{L^\infty(\Omega)} = \dis \inf_{\varphi} \sup_{\Omega} |H|^2\), the proof is completed.
	\end{proof}
	For the hyperbolic space \(\mathbb{H}^n(-1)\), Cheng and Yang do not rely on Nash's theorem; instead, they construct an appropriate trial function in \cite{ChengYang11} to derive a universal inequality for eigenvalues of the eigenvalue problem \eqref{DL}.
	\begin{thm}
		For the eigenvalue problem of the Dirichlet Laplacian \eqref{DL} on a domain in the hyperbolic space \(\mathbb{H}^n(-1)\), eigenvalues \(\lambda_k\) satisfy
		\begin{equation}\label{ineq:ChengYang-hyper}
			\sum_{i = 1}^{k} (\lambda_{k + 1} - \lambda_i)^2 \leq 4 \sum_{i = 1}^{k} (\lambda_{k + 1} - \lambda_i) \left(\lambda_i - \frac{(n - 1)^2}{4}\right).	
		\end{equation}
	\end{thm}
	According to the recursion formula in the theorem \ref{thm:ChengYang07} of Cheng and Yang, we know that 
	\begin{thm}\label{thm:ChengYangsub}
		Let \(\Omega\) be a bounded domain in an \(n\)-dimensional complete Riemannian manifold \(M^n\). There exists a constant \(H_0^2\), which depends only on \(M\) and \(\Omega\), such that eigenvalues \(\lambda_k\)'s of the eigenvalue problem \eqref{DL} satisfy
		\begin{equation}\label{ineq:ChengYang-LY1}
			\frac{1}{k} \sum_{i = 1}^{k} \lambda_i + \frac{n^2}{4} H_0^2 \geq \frac{n}{\sqrt{(n + 2)(n+4)}} C_n \frac{k^{\frac{2}{n}}}{\abs{\Omega}^{2/n}}, \text{ for } k = 1, 2, \dots,
		\end{equation}
		where 
		\begin{equation*}
			C_n = \frac{4 \pi^2}{\omega_{n}^{2/n}}.
		\end{equation*}
	\end{thm}
	Furthermore, Cheng and Yang \cite{ChengYang11} propose the following conjecture.
	\begin{conjecture}[{\cite{ChengYang11}}]\label{conj:ChengYang11}
		Let \(\Omega\) be a bounded domain in an \(n\)-dimensional complete Riemannian manifold \(M^n\). Then, there exists a constant \(c(M, \Omega)\), which depends only on \(M\) and \(\Omega\), such that eigenvalues \(\lambda_{k}\) satisfy
		\begin{equation}\label{ineq:ChengYangcon1}
			\frac{1}{k} \sum_{i=1}^{k} \lambda_{i} + c(M, \Omega) \geq \frac{n}{n+2} C_n \frac{k^{\frac{2}{n}}}{\abs{\Omega}^{2/n}}, \quad \text{for } k = 1, 2, \dots,
		\end{equation}
		\begin{equation}\label{ineq:ChengYangcon2}
			\lambda_{k} + c(M, \Omega) \geq C_n \frac{k^{\frac{2}{n}}}{\abs{\Omega}^{2/n}}, \quad \text{for } k = 1, 2, \dots.	
		\end{equation}
	\end{conjecture}
	\begin{thm}\label{ChenCheng24}
		Let \(\Omega\) be a bounded domain in an \(n\)-dimensional complete Riemannian manifold \(M^n\). Then, there exists a constant \(H_0\), which depends only on \(M\) and \(\Omega\), such that eigenvalues \(\lambda_{i}\) satisfy
		\begin{equation}\label{ineq:ChenCheng24}
			Z_H(t) \leq \dfrac{|\Omega|}{(4 \pi t)^{n/2}},		
		\end{equation}
		where
		\begin{equation*}
			Z_H(t) = \dis \sum_{i=1}^{\infty} \exp \br{-\br{\lambda_i + \dfrac{n^2}{4} H_0^2} t}.
		\end{equation*}
	\end{thm}
	Since the projective spaces admit the first standard embedding into Euclidean space (See Section \ref{sec:pre}), by the Cheng-Yang recursion formula \eqref{ineq:ChengYang-LY1}, we have
	\begin{cor}\label{cor:PS}
		Let \(M\) be an \(n\)-dimensional submanifold in the projective space \(\mathbb{F}P^m\) with the mean curvature vector fields \(H\) and let \(\Omega\) be a bounded domain in \(M^n\). Then eigenvalues \(\lambda_k\)'s of the eigenvalue problem \eqref{DL} satisfy	
		\begin{equation}\label{ineq:ChenCheng24PS}
			\frac{1}{k} \sum_{i=1}^{k} \lambda_{i} + \dfrac{n^2}{4} \br{H_0^2 + \frac{2(n + d(\mathbb{F}))}{n}} \geq \frac{n}{\sqrt{(n + 2)(n+4)}} C_n \frac{k^{\frac{2}{n}}}{\abs{\Omega}^{2/n}}, \quad \text{for } k = 1, 2, \dots,
		\end{equation}
		where \(d(\mathbb{F})\) is defined by \eqref{defn:dim} in Section \ref{sec:pre} and \(H_0^2 = \|H\|^2_{L^\infty(\Omega)} = \sup_{\Omega} |H|^2\). In particular, if \(M\) is an \(n\)-dimensional minimal submanifold in the projective space \(\mathbb{F}P^m\), then eigenvalues satisfy
		\begin{equation}\label{ineq:ChenCheng24PS0}
			\frac{1}{k} \sum_{i=1}^{k} \lambda_{i} + \dfrac{n (n + d(\mathbb{F}))}{2} \geq \frac{n}{\sqrt{(n + 2)(n+4)}} C_n \frac{k^{\frac{2}{n}}}{\abs{\Omega}^{2/n}}, \quad \text{for } k = 1, 2, \dots,
		\end{equation}
	\end{cor}
	\begin{rem}
		For compact submanifolds with boundaries immersed in Euclidean spaces, spheres, and projective spaces, we consider the closed eigenvalue problem. The estimates for eigenvalue inequalities presented in the aforementioned theorems remain valid. The methods of proof and computations require no modification, except that eigenvalues now begin with \(\lambda_1 = 0\). In particular, for the minimal submanifold \(M^n\) without boundary immersed into \(\overline{M}\), where \(\overline{M}\) is the unit sphere \(\mathbb{S}^{N}(1)\) or one of projective spaces \(\mathbb{F}P^m\), the closed eigenvalues of \(M^n\) obey
		\begin{equation*}
			\frac{1}{k} \sum_{i=1}^k \lambda_i + \frac{c(n)}{4} \geq \frac{n}{\sqrt{(n + 2)(n+4)}} C_n \frac{k^{\frac{2}{n}}}{\abs{M}^{2/n}}, \quad \text{for } k = 1, 2, \dots,	
		\end{equation*}
		where \(c(n)\) depends only on dimension \(n\) defined by
		\begin{equation*}
			c(n) = \begin{cases}
				n^2, & \text{for } \overline{M} \text{ is } \mathbb{S}^{N}(1), \\
				2n(n + d(\mathbb{F})), & \text{for } \overline{M} \text{ is } \mathbb{F}P^m.
			\end{cases}
		\end{equation*}
		Since the projective spaces \(\mathbb{F}P^m\) may be minimally embedded into the sphere \(\mathbb{S}^{(m+1)d(\mathbb{F})-1} \br{\sqrt{\dfrac{m}{2(m+1)}}}\) (See \cite{ChenBY84,Tai68}), the closed eigenvalues of \(\mathbb{F}P^m\) satisfy
		\begin{equation*}
			\frac{1}{k} \sum_{i=1}^k \lambda_i + \frac{m(m+1)d^2(\mathbb{F})}{2} \geq \frac{m d(\mathbb{F})}{\sqrt{(m d(\mathbb{F}) + 2)(m d(\mathbb{F}) + 4)}} C_{m d(\mathbb{F})} \frac{k^{\frac{2}{m d(\mathbb{F})}}}{\abs{M}^{\frac{2}{m d(\mathbb{F})}}}, \quad \text{for } k = 1, 2, \dots.
		\end{equation*}
	\end{rem}
	
	\section{Proof of Theorem \ref{ChenCheng24} and Corollary \ref{cor:PS}}\label{sec:pre}
	
	In order to prove Theorem \ref{ChenCheng24} and Corollary \ref{cor:PS}, we require the following Tauberian theorem and the first standard embedding of projective spaces into Euclidean spaces.
	\subsection{Karamata's Tauberian Theorem}
	\begin{thm}\label{them:WeylLaw1}\cite{LMP23}
		Let \(M\) be an \(n\)-dimensional smooth compact Riemannian manifold. If \(\partial M \neq \varnothing\), assume that either the Dirichlet or the Neumann boundary conditions are imposed on the boundary. Then the eigenvalue counting function for \(M\) has the asymptotics
		\begin{equation}\label{Weyl-asymptotic}
			\mathcal{N}(\lambda) = \Lon \mathrm{Vol}(M) \lambda^{\frac{n}{2}} + O\left(\lambda^{\frac{n - 1}{2}}\right),	
		\end{equation}
		where \(\Lon\) is defined by 
		\begin{equation}\label{defn:classical-constant-0}
			\Lc_{0,n}^{cl} := \omega_n / (2\pi)^n 
		\end{equation}
		is called the classical constant and \(\omega_n\) is the volume of the \(n\)-ball, \(\omega_n = \pi^{n/2} / \Gamma(1 + n/2)\).
	\end{thm}
	A fundamental property for \(\rho, \delta > 0\), often referred to as the Riesz iteration or alternatively as the Aizenman-Lieb procedure \cite{AL78}, states that
	\begin{equation}\label{RieszIteration}
		R_{\rho + \delta}(\lambda) = \dfrac{\Gamma(\rho + \delta + 1)}{\Gamma(\rho + 1) \Gamma(\delta)} \int_{0}^{\infty} \left(\lambda - t\right)_{+}^{\delta - 1} R_{\rho}(t) \, dt.
	\end{equation}
	It is noteworthy that applying the Riesz iteration \eqref{RieszIteration} to (\ref{Weyl-asymptotic}) directly yields the following result
	\begin{equation}\label{riesz-weyl}
		R_{\rho}(z) \sim \Lrn |\Omega| z^{\rho + n/2}, \quad \text{as } z \to \infty, 
	\end{equation}
	where the classical constant is given by \eqref{defn:classical-constant}.
	This formula allows us to estimate the counting function \(N(\lambda) = \# \{n : \lambda_n < \lambda\}\) with the help of the kernel \(Z(t)\). For this, we will make use of the following Tauberian theorem due to Karamata \cite{Karamata} (See also \cite[Theorem 1.1]{ANPT09}).
	\begin{thm}\label{thm:Karamata}
		Let \((\lambda_n)_{n \in \mathbb{N}}\) be a sequence of positive real numbers such that the series \(\dis \sum_{n \in \mathbb{N}} e^{-\lambda_n t}\) converges for every \(t > 0\). Then for \(r > 0\) and \(a \in \mathbb{R}\), the following are equivalent
		\begin{enumerate}[(1)]
			\item \(\dis \lim_{t \to 0} t^r \sum_{n \in \mathbb{N}} e^{-\lambda_n t} = a\),
			\item \(\dis \lim_{\lambda \to \infty} \lambda^{-r} N(\lambda) = \frac{a}{\Gamma(r + 1)}\).
		\end{enumerate}
		Here \(N\) denotes the counting function \(N(\lambda) = \# \{\lambda_n \leq \lambda\}\), and \(\Gamma(r) = \dis \int_{0}^{\infty} x^{r - 1} e^{-x} \, dx\) is the usual Gamma function. 	
	\end{thm}
	
	\subsection{Submanifolds in Projective Spaces}\label{sub:PS}
	Let \(\mathbb{F}\) denote the field \(\mathbb{R}\) of real numbers, the field \(\mathbb{C}\) of complex numbers, or the field \(\mathbb{Q}\) of quaternions. Let \(\mathbb{F}P^m\) denote the \(m\)-dimensional projective space over \(\mathbb{F}\). The projective space \(\mathbb{F}P^m\) is endowed with a standard Riemannian metric whose sectional curvature is either constant and equal to \(1\) (for \(\mathbb{F} = \mathbb{R}\)) or pinched between \(1\) and \(4\) (for \(\mathbb{F} = \mathbb{C}\) or \(\mathbb{Q}\)). It is well-known that projective spaces admit a first standard embedding into Euclidean spaces (see, e.g., \cite{ChenBY84, SoufiHaIlias09, Tai68}). Let \(\rho: \mathbb{F}P^m \to H_{m+1}(\mathbb{F})\) denote the first standard embedding of projective spaces into Euclidean spaces, where
	\[
	H_{m+1}(\mathbb{F}) = \{ A \in M_{m+1}(\mathbb{F}) \mid A = A^t \},
	\]
	and \(M_{m+1}(\mathbb{F})\) denotes the space of \((m+1) \times (m+1)\) matrices over \(\mathbb{F}\). For convenience, we introduce the integer
	\begin{equation}\label{defn:dim}
		d(\mathbb{F}) = \dim_{\mathbb{R}} \mathbb{F} =
		\begin{cases}
			1, & \text{if } \mathbb{F} = \mathbb{R}, \\
			2, & \text{if } \mathbb{F} = \mathbb{C}, \\
			4, & \text{if } \mathbb{F} = \mathbb{Q}.
		\end{cases}	
	\end{equation}
	\begin{prop}\label{prop:proj}\cite[Proposition 2.4]{Chen09}
		Let \(f: M^n \to \mathbb{F}P^m\) be an isometric immersion, and let \(H\) and \(H'\) denote the mean curvature vector fields of the immersions \(f\) and \(\rho \circ f\), respectively. Then,
		\[
		\abs{H'}^2 = \abs{H}^2 + \frac{4(n + 2)}{3n} + \frac{2}{3n^2} \sum_{i \neq j} \tilde{K}(e_i, e_j),
		\]
		where \(\{e_i\}_{i=1}^n\) is a local orthonormal basis of \(TM\), and \(\tilde{K}\) is the sectional curvature of \(\mathbb{F}P^m\), which can be expressed as
		\[
		\tilde{K}(e_i, e_j) =
		\begin{cases}
			1, & \text{if } \mathbb{F} = \mathbb{R}, \\
			1 + 3 (e_i \cdot J e_j)^2, & \text{where } J \text{ is the complex structure of } \mathbb{C}P^m, \text{ if } \mathbb{F} = \mathbb{C}, \\
			1 + 3 \dis \sum_{r=1}^3 (e_i \cdot J_r e_j)^2, & \text{where } J_r \text{ is the quaternionic structure of } \mathbb{Q}P^m, \text{ if } \mathbb{F} = \mathbb{Q}.
		\end{cases}
		\]
		In particular, it follows that
		\[
		|H'|^2 =
		\begin{cases}
			|H|^2 + \dfrac{2(n + 1)}{n}, & \text{for } \mathbb{R}P^m, \\
			|H|^2 + \dfrac{2(n + 1)}{n} + \dfrac{2}{n^2} \dis \sum_{i,j=1}^n (e_i \cdot J e_j)^2 \leq |H|^2 + \dfrac{2(n + 2)}{n}, & \text{for } \mathbb{C}P^m, \\
			|H|^2 + \dfrac{2(n + 1)}{n} + \dfrac{2}{n^2} \dis \sum_{i,j=1}^n \sum_{r=1}^3 (e_i \cdot J_r e_j)^2 \leq |H|^2 + \dfrac{2(n + 4)}{n}, & \text{for } \mathbb{Q}P^m.
		\end{cases}
		\]
		That is,
		\[
		|H'|^2 \leq |H|^2 + \frac{2(n + d(\mathbb{F}))}{n},
		\]
		where equality holds if and only if \(M\) is a complex submanifold of \(\mathbb{C}P^m\) (for the case of \(\mathbb{C}P^m\)), or \(n \equiv 0 \pmod{4}\) and \(M\) is a quaternionic submanifold of \(\mathbb{Q}P^m\) (for the case of \(\mathbb{Q}P^m\)).
	\end{prop}
	
	\subsection{Harrell-Stubbe Type Estimates}
	To prove Theorem \ref{ChenCheng24}, we require the following Harrell-Stubbe type estimates.
	\begin{thm}\label{CCRiesz}
		Let \(\Omega\) be a bounded domain in an \(n\)-dimensional complete Riemannian manifold \(M^n\). For the eigenvalue problem of the Dirichlet Laplacian (\ref{DL}), there exists a constant \(H_0^2\), which depends only on \(M\) and \(\Omega\), such that
		\begin{equation}\label{ineq:Riesz1}
			R_2(z) \leq \frac{4}{n} \sum_{k \in \mathbb{N}} \br{z - \mu_i}_+ \mu_i,
		\end{equation}	
		where \(\mu_i\) is defined by 
		\begin{equation}\label{defn:mu}
			\mu_i = \lambda_i + \frac{n^2}{4} H_0^2.
		\end{equation}
		Furthermore, for \(\rho \geq 2\) and \(z \geq \mu_1\),
		\begin{equation}\label{Riesz2}
			R_\rho(z) \leq \frac{\rho}{\rho + \frac{n}{2}} z R_{\rho-1}(z).
		\end{equation}	
	\end{thm}
	\begin{rem}
		If the Berezin-Li-Yau \eqref{Riesz2} holds for \(\rho \geq 1\), then the inequalities \eqref{ineq:ChengYangcon1} in Conjecture \ref{conj:ChengYang11} are also satisfied.
	\end{rem}
	\begin{proof}[Proof of Theorem \ref{CCRiesz}]
		The inequality \eqref{ineq:Riesz1} can be established using a similar argument to that employed in the proof of Proposition 1 in \cite{ChengYang06} and Theorem 1.1 in \cite{ChenCheng08}. Assume that \(u_i\) is an orthonormal eigenfunction corresponding to the \(i\)-th eigenvalue \(\lambda_i\), i.e., \(u_i\) satisfies
		\begin{equation}\label{defn:Chencheng1}
			\begin{cases}
				\Delta u_i = -\lambda_i u_i, & \text{in } M, \\
				\left.u_i\right|_{\partial M} = 0, \\
				\dis \int_\Omega u_i u_j \, d\mu = \delta_{ij}.
			\end{cases}	
		\end{equation}
		Define \(\varphi_i, a^\alpha_{ij}\), and \(b^\alpha_{ij}\), for \(i, j = 1, \dots, k\), by
		\begin{equation}\label{defn:Chencheng2}
			\begin{cases}
				a^\alpha_{ij} = \dis \int_M y^\alpha u_i u_j \, d\mu, \\
				\varphi_i^\alpha = y^\alpha u_i - \dis \sum_{j = 1}^{k} a^\alpha_{ij} u_j, \\
				b^\alpha_{ij} = \dis \int_\Omega u_j \left( \nabla u_i \cdot \nabla y^\alpha + \frac{1}{2} u_i \Delta y^\alpha \right) \, d\mu.
			\end{cases}	
		\end{equation}
		Then we can deduce from \eqref{defn:Chencheng1} and \eqref{defn:Chencheng2}
		\begin{equation*}
			\begin{aligned}
				& a^\alpha_{ij} = a^\alpha_{ji}, \quad 2 b^\alpha_{ij} = -2 b^\alpha_{ji} = \br{\lambda_{i} - \lambda_j} a^\alpha_{ij}, \\
				& \int_\Omega \varphi^\alpha_i u_j \, d\mu = 0, \text{ for } j = 1, 2, \dots, k.	
			\end{aligned}
		\end{equation*}
		For \(\lambda_{k} \leq z < \lambda_{k+1}\), from the Rayleigh-Ritz inequality, we have
		\begin{equation*}
			z < \lambda_{k + 1} \leq \dfrac{\dis \int_\Omega \abs{\nabla \varphi^\alpha_i}^2 \, d\mu}{\dis \int_\Omega \br{\varphi^\alpha_i}^2 \, d\mu}.	
		\end{equation*}
		Taking a similar argument as in \cite{ChengYang06} (see also \cite{AshbaughHermi04, AshbaughHermi07, ChengYang07}), we can infer
		\begin{equation}\label{ineq:Chencheng1}
			\sum_{i = 1}^{k} (z - \lambda_i)^2 \|u_i \nabla y^\alpha\|^2 \leq \sum_{i = 1}^{k} (z - \lambda_i) \|2 \nabla y^\alpha \cdot \nabla u_i + u_i \Delta y^\alpha\|^2.
		\end{equation}
		From Lemma \ref{lem:ChenCheng08}, a direct calculation yields
		\begin{equation}\label{ineq:Chencheng2}
			\begin{aligned}
				\sum_{\alpha=1}^N \|u_i \nabla y^\alpha\|^2 =& n, \\
				\sum_{\alpha=1}^N \norm{2 \nabla y^\alpha \cdot \nabla u_i + u_i \Delta y^\alpha}^2 =& 4 \lambda_i + n^2 \int_{\Omega} \abs{H}^2 u_i^2 \, d\mu \\
				\leq & 4 \br{\lambda_i + \dfrac{n^2}{4} H_0^2}.
			\end{aligned}	
		\end{equation}
		Putting \eqref{ineq:Chencheng2} into \eqref{ineq:Chencheng1} yields, for \(\lambda_k \leq z < \lambda_{k+1}\),
		\begin{equation}\label{ineq:AshHermi5}
			\sum_{i=1}^{k} \br{z - \lambda_i}^2 \leq \frac{4}{n} \sum_{i=1}^{k} \br{z - \lambda_i} \br{\lambda_i + \dfrac{n^2}{4} H_0^2}.
		\end{equation}
		which implies \eqref{ineq:Riesz1}. 
		From \eqref{ineq:Riesz1}, we infer
		\begin{equation*}
			\begin{aligned}
				R_2(z) \leq & \frac{4}{n} \sum_{k \in \mathbb{N}} \br{z - \mu_i}_+ \mu_i \\
				= & \frac{4}{n} \sum_{k \in \mathbb{N}} \br{z - \mu_i}_+ \br{z - \br{z - \mu_i}_+} \\
				= & \frac{4}{n} z R_1(z) - \frac{4}{n} R_2(z),
			\end{aligned}	
		\end{equation*}
		which implies 
		\begin{equation}\label{ineq:AshHermi7}
			\br{1 + \frac{n}{4}} R_2(z) \leq z R_1(z).	
		\end{equation}
		The inequality \eqref{Riesz2} can be deduced from \eqref{ineq:AshHermi7} by using Riesz iteration, which completes the proof of Theorem \ref{CCRiesz}.
	\end{proof}
	
	\subsection{Proof of Theorem \ref{ChenCheng24} and Corollary \ref{cor:PS}}
	\begin{proof}[Proof of Theorem \ref{ChenCheng24}]
		From the inequality \eqref{Riesz2} in Theorem \ref{CCRiesz}, the function
		\begin{equation}\label{ineq:monod}
			z \mapsto \frac{R_{\rho}(z)}{z^{\rho + \frac{n}{2}}}
		\end{equation}
		is a nondecreasing function of \(z\), for \(\rho \geq 2\). According to the asymptotic formula \eqref{riesz-weyl} and Theorem \ref{thm:Karamata}, we get the Berezin-Li-Yau inequality 
		\begin{equation}\label{ineq:BLY}
			R_{\rho}(z) \leq L_{\rho,n}^{cl} |\Omega| z^{\rho + n/2},
		\end{equation}
		where \(\Lrn\) is defined in \eqref{defn:classical-constant}. The Laplace transform yields
		\begin{equation}\label{equ:Lap}
			\mathcal{L} \bigl( (z - \mu_k)_+^\rho \bigr) = \frac{\Gamma(\rho + 1) e^{-\mu_k t}}{t^{\rho + 1}}.
		\end{equation}
		Combining the inequality \eqref{ineq:BLY} with \eqref{equ:Lap} yields
		\[
		\frac{\Gamma(\rho + 1)}{t^{\rho + 1}} Z_H(t) \leq L_{\rho,n}^{cl} |\Omega| \frac{\Gamma \bigl( \rho + 1 + \frac{n}{2} \bigr)}{t^{\rho + 1 + \frac{n}{2}}},
		\]
		where \(Z_H(t)\) is defined by	
		\[
		Z_H(t) = \dis \sum_{i=1}^{\infty} \exp \br{-\mu_i t}.
		\]
		From the definition of \(\Lrn\) in \eqref{defn:classical-constant}, we obtain the desired result 
		\begin{equation}\label{ineq:KacCC}
			Z_H(t) \leq \dfrac{|\Omega|}{(4 \pi t)^{n/2}}.
		\end{equation}
		This completes the proof of Theorem \ref{ChenCheng24}.
	\end{proof} 
	\begin{proof}[Proof of Corollary \ref{cor:PS}]
		From Proposition \ref{prop:proj}, the projective spaces \(\mathbb{F}P^m\), which admit the first standard embedding into the Euclidean space \(H(m+1, \mathbb{F})\), allow an isometrically immersed submanifold \(M^n\) in \(\mathbb{F}P^m\) to be regarded as a submanifold in Euclidean space \(H(m+1, \mathbb{F})\). The mean curvature vector field \(H'\) of \(M^n\) in Euclidean space satisfies
		\[
		|H'|^2 \leq |H|^2 + \frac{2(n + d(\mathbb{F}))}{n},
		\]
		where \(H\) denotes the mean curvature vector field of \(M^n\) in \(\mathbb{F}P^m\). Consequently, Corollary \ref{cor:PS} follows the analogous argument of the proof of Theorem \ref{thm:ChenCheng08submfd}.	
	\end{proof}
	
	\textbf{Acknowledgement.} This work was partly supported by MEXT Promotion of Distinctive Joint Research Center Program JPMXP0723833165 and Osaka Metropolitan University Strategic Research Promotion Project (Development of International Research Hubs).

	\begin{flushleft}
		Daguang Chen \\
		E-mail: dgchen@mail.tsinghua.edu.cn \\
		Department of Mathematical Sciences \\
		Tsinghua University, Beijing, 100084, China 	
	\end{flushleft}
	\vspace{0.5cm}
	\begin{flushleft}
		Qing-Ming Cheng \\
		E-mail: qingmingcheng@yahoo.com, chengqingming@cqut.edu.cn \\
		Mathematical Science Research Center \\
		Chongqing University of Technology, Chongqing, 400054, China \\
		Osaka Center Advanced Mathematical Institute \\
		Osaka Metropolitan University, Osaka 558-8585, Japan 
	\end{flushleft}
	
\end{document}